\documentclass[12pt]{article}
\usepackage{amsmath}
\usepackage{amsthm}
\usepackage{amssymb,graphics,fancybox}
\usepackage[final]{graphicx}

\def\PA{Painlev\'e }

\newtheorem{theorem}{Theorem}[section]
\newtheorem{lemma}[theorem]{Lemma}

\theoremstyle{definition}

\theoremstyle{remark}

\title{R.~Fuchs' problem of the Painlev\'e equations from the first to the fifth}
\author{Yousuke Ohyama \\
\small  Graduate School of Information Science and Technology,\\
\small Osaka University\\ {}\\
Shoji Okumura\\
\small  Graduate School of Science, Osaka University\\
}
\begin{document}
\maketitle
%% \tableofcontents

\section{Introduction}
About one hundred  years have passed since R.~Fuchs showed that 
the sixth \PA equation is represented as a monodromy preserving deformation \cite{Fuchs1}. 
Garnier showed that every type of the \PA equations is also 
represented as a monodromy preserving deformation \cite{Garnier}. 
Although the monodromy data is preserved by the deformations,  
it is difficult to calculate the monodromy data in general. 

Riemann calculated the monodromy group for the Euler-Gauss hypergeometric equation.
Moreover,   monodromy groups of the Euler-Gauss hypergeometric equation 
are the polyhedron  groups when it has algebraic solutions \cite{Schw}. 
Schwarz' fifteen algebraic solutions are obtained from  simple hypergeometric 
equations by rational transformations of the independent variable. 
We will study an analogue of Schwarz' solutions in the \PA case.

In this paper we call a linear equation  whose isomonodromic 
deformation gives a \PA equation as {\it the \PA type}. 
For  the linear equation of the \PA type, we can calculate the monodromy data 
explicitly if we substitute particular solutions of the \PA equations.
Historically, the first example of such solutions  is 
Picard's solutions of the sixth \PA  equation \cite{Pi}.
R.~Fuchs calculated the monodromy group of the linear equation 
corresponding to Picard's solutions   \cite{Fuchs:11}. 
It seems that R.~Fuchs' paper \cite{Fuchs:11}, whose title is  the same  as 
the famous paper \cite{Fuchs1}, has been  forgotten for long years.
Recently Mazzocco found again his result independently but 
his paper was not referred in her paper  \cite{Ma1}. 
In \cite{Fuchs:11}, R.~Fuchs proposed the following problem: 
\begin{quote}
{\bf R.~Fuchs' Problem}\quad 
Let $y(t)$ be an algebraic solution $y(t)$ of a \PA equation. 
Find a suitable transformation $x=x(z,t)$ such that the corresponding 
linear differential equation 
$$\frac{d^2 v}{dz^2}=Q(t,y(t),z)v$$
is changed  to the form without the deformation parameter $t$:
$$\frac{d^2 u}{dx^2}=\tilde{Q}(x)u.$$
Here $v= \sqrt{ {dz}/{dx}}\ u.$ (See the lemma \ref{keylemma}.)
\end{quote}
Picard's solution is algebraic if it corresponds to rational points of 
elliptic curves. For three, four and six divided points, 
the genus of algebraic Picard's solutions is zero. 
 R.~Fuchs showed that for algebraic Picard's solutions whose genus are zero  
the corresponding linear equations are reduced the Euler-Gauss hypergeometric 
equations by suitable rational transformations.
 
If R.~Fuchs' problem is true, algebraic solutions of the \PA equations can be considered as
a kind of a generalization of Schwarz' solutions.  Schwarz' solutions are 
constructed by rational transformations which change hypergeometric equations 
to hypergeometric equations. Algebraic solutions of the \PA equations are 
constructed by rational transformations which change hypergeometric equations 
to linear equations of the \PA type and we can calculate the monodromy data of 
the linear equation explicitly.

Recently, Kitaev, jointing with Andreev,   constructed many algebraic solutions of the 
sixth \PA equation, which include  known ones and new ones, by rational transformations 
 of the hypergeometric equations \cite{AK}, \cite{Kitaev:05}.  
At least now, we do not know whether R.~Fuchs' problem is true or not for the 
sixth \PA equation. We do not have negative example for R.~Fuchs' problem. 
Kitaev's transformation is a generalization of classical work by Goursat 
\cite{Goursat}. 
Goursat found  many rational transformations which keep the hypergeometric equations. 
His list is incomplete and Vid\=unas made a complete list of rational transformations \cite{Vidunas}.

In this paper we will study a confluent version of \cite{AK} and 
we show that  R,~Fuchs' problem is true for the \PA equation from the first 
to the fifth.  
We will classify all  rational  transformations which change 
the confluent hypergeometric equations to linear equations  of  
the \PA type from the first to the fifth.  Compared with the sixth \PA equation,
we obtain less transformations since the Poincar\'e rank of irregular singularities
of the linear equations 
of the \PA type is up to three. We have up to sixth order rational 
functions which change the confluent hypergeometric equations to linear equations of 
the \PA type. 

We show such rational transformations correspond  to almost all  of algebraic solutions 
of the \PA equations from the first to fifth up to   the B\"acklund transformations.
The cases of the degenerate fifth \PA equation 
and the Laguerre type solution of the  fifth \PA equation are  exceptional. 
We need exponential type transformations since the monodromy data is decomposable. 

The parabolic cylinder equation (the Weber equation),
the Bessel equation and the Airy equation are reduced to the confluent hypergeometric 
equations by rational transformations of the independent variable (see \eqref{bessel},
\eqref{weber} and \eqref{airy}). These classical formula can be considered 
as confluent version of the  Goursat transformations. Our result is 
an analogue of such classical formula  in the \PA analysis.

Moreover, we obtain four non-algebraic solutions 
by rational transformations from the confluent hypergeometric equations.
They are called as symmetric solutions of the \PA equations 
which are fixed points of simple transformations of the \PA equations. 
The symmetric solutions are not classical solutions in Umemura's sense, 
but it is a kind of generalization of classical solutions.

In the section \ref{PA}, we will review the \PA equations and their special solutions. 
In the section \ref{chg}, we will review confluent hypergeometric equations since 
we treat degenerate confluent hypergeometric equations. 
In the section \ref{iso}, we will give a list of linear equations which correspond 
to the \PA equations. In the section  \ref{rational}, we will show that 
R.~Fuchs' problem is true for the \PA equations from the first to the fifth.

The authors give   thanks to Professor Kazuo Okamoto, Professor Alexander Kitaev and Mr. Kazuo Kaneko for 
fruitful discussions. 

\section{Remarks on the \PA equations}\label{PA}
In order to fix the notation, we will review the \PA equations.

\subsection{Remarks on P2 and P34}
The thirty-fourth \PA equation P34($a$)
$$y^{\prime\prime}=\frac{(y^{\prime })^2}{2y}+2y^2 -ty -\frac{a}{2y} $$
appeared in Gambier's list of second order nonlinear equations without 
movable singularities \cite{Gambier}. It is known that P34 is equivalent to
the second \PA equation P2($\alpha$)
$$y^{\prime\prime}= 2y^3 + t y +\alpha.$$
P2($\alpha$) can be written in the Hamiltonian form 
\begin{equation}\label{2:ha2}
{\cal H}_{II}: \
\left\{\begin{array}{l}
   q^{\prime } = -q^2 +p -\frac t2, \\
   p^{\prime } = 2pq+ \left(\alpha+\frac12\right), \\
\end{array}
\right.
\end{equation}
with the Hamiltonian
 $$ H_{II}= {1 \over 2}p^2- \left(q^2+{t \over 2} \right) p-\left(\alpha+\frac12\right) q.$$
If we eliminate $p$ from \eqref{2:ha2}, we obtain P2($\alpha$). 
If we eliminate $q$ from \eqref{2:ha2}, we obtain the thirty-fourth 
\PA equation P34($(\alpha+1/2)^2$).
Therefore P2($\alpha$) and P34($(\alpha+1/2)^2$) are equivalent as nonlinear 
differential equations, but we will distinguish these two equations from the isomonodromic viewpoint. 

\subsection{Remarks on P3 and P5}

In the following we distinguish the three types of 
the third \PA equation P3$(\alpha, \beta, \gamma, \delta)$:
\begin{eqnarray*}
 y^{\prime\prime}&=& \frac{1}{ y}{y^{\prime}}^2-\frac{y^{\prime}}{ t}+
\frac{\alpha y^2 + \beta} {t}+\gamma y^3 + \frac{\delta }{ y},
\end{eqnarray*}
because they have different forms of isomonodromic deformations:
\begin{itemize}
	\item $D_8^{(1)}$ if $\alpha\not=0, \beta\not=0$, $\gamma=0, \delta=0$,
	\item $D_7^{(1)}$ if $\delta=0, \beta\not=0$ or $\gamma=0,  \alpha\not=0$,
	\item $D_6^{(1)}$ if $\gamma\delta\not=0$.
\end{itemize}
In the case $\beta=0, \delta=0$ (or $\alpha=0, \gamma=0$), the third \PA equation 
is a quadrature, and we exclude this case from the \PA family.  $D_j^{(1)} (j=6,7,8)$ mean 
the affine Dynkin diagrams corresponding to Okamoto's initial value spaces.  
In this paper we omit the upper index $(1)$ for simplicity and denote $D_6, D_7, D_8$.
For details, see \cite{KOOS}.
 By   suitable scale transformations
$t \to ct, y \to dt$, we can fix $\gamma=4, \delta=-4$ for $D_6$. 
and $\gamma=2$ for $D_7$. 
We will use another form of the third \PA equation P3${}^\prime(\alpha, \beta, \gamma, \delta)$ 
\begin{eqnarray*}
q^{\prime\prime}&=& 
\frac{1}{q}{q^{\prime}}^2-\frac{q^{\prime}}{x}+
\frac{q^2(\alpha+\beta)} {4x^2}+\frac{\gamma}{4x} + \frac{\delta }{4q},
\end{eqnarray*}
since P3${}^\prime$ is more sympathetic to isomonodromic deformations than P3. 
We can change P3 to P3${}^\prime$ by $x=t^2, ty = q.$
\par\bigskip

For the fifth equations P5$(\alpha, \beta, \gamma, \delta)$ 
\begin{eqnarray*}
y^{\prime\prime}&=& \left(\frac{1}{2y}+\frac{1}{y-1} \right){y^{\prime}}^2-\frac{1}{t}
{y^{\prime}}+\frac{(y-1)^2}{t^2}\left(\alpha y +\frac{\beta}{y}\right)\\
&{}& \hskip 2cm +\gamma {y\over t}+\delta {y(y+1)\over y-1},
\end{eqnarray*}
we assume that $\delta\not=0$. 
When $\delta=0, \gamma\not=0$, the fifth equation is equivalent to the third equation of 
the $D_6$ type \cite{Gromak}. We denote deg-P5$(\alpha, \beta, \gamma, 0)$ if $\delta=0$. 
By a suitable scale transformation 
$t \to ct$, we can fix $\delta=-1/2$ for P5 and  $\gamma=-2$ for deg-P5. 
Let $q$ is a solution of  P3${}^\prime (4(\alpha_1-\beta_1),-4(\alpha_1+\beta_1-1), 4, \ -4)$.
Then
$$ y =\frac{tq'-q^2-(\alpha_1+\beta_1)q-t}{tq'+q^2-(\alpha_1+\beta_1)q-t}$$
is a solution of  deg-P5$( {\alpha_1^2}/2, -{\beta_1^2}/2, -2,0)$.

P5$(\alpha, \beta, 0, 0)$ is quadrature and we exclude this case from the \PA family.

\subsection{Special solutions from P1 to P5}
We will review special solutions of the  \PA equations. 
Although generic solutions of the  \PA equations are transcendental, 
there exist some special solutions which reduce to known functions.
Umemura defined a class of `known functions', which are
called {\it classical solutions} \cite{U:stras}. He also gave a method 
how to classify classical solutions of the \PA equations. 

There exist two types of classical solutions of \PA equations, 
one is algebraic and the second is the Riccati type. 
Up to now, all of classical solutions  are classified 
except algebraic  solutions of the sixth \PA equations. 
We will list up all of classical solutions of P1 to P5. 

\begin{theorem}\label{classical}
1) All solutions of P1 are transcendental.

\par\medskip\noindent
2) P2($0$) has a rational solution $y=0$. 
P2($-1/2$) has a Riccati type solution $y=-u'/u$. 
Here $u$ is any solution of the Airy equation $ u^{\prime\prime}+tu/2=0$.

\par\medskip\noindent
2) P34($(\alpha+1/2)^2$)  is equivalent to P2($\alpha$).  
P34($1/4$) has a rational solution $y=t/2$. P34($1$) has the Riccati type solutions. 

\par\medskip\noindent
4) P4($0,-2/9$) has a rational solution $y=-2t/3$.  
P4($1-s, -2s^2$) has a Riccati type solution $y=-u'/u$. 
Here $u$ is any solution of the Hermite-Weber equation 
$ u^{\prime\prime}+2t u^{\prime}+2s u  =0$. If $s=1$, 
P4($0,-2$) has a rational solution $y=-2t$, which reduces to the 
Hermite polynomial.

\par\medskip\noindent
5) P3${}^{\prime}(D_6)$($a,-a, 4,-4$)  has an algebraic solution $y=-\sqrt{t}$.
P3$^{\prime}(D_6)$($4h, 4(h+1),4,-4$) has a Riccati type solution $y=u'/u$.
Here $u$ is any solution  of
$ tu^{\prime\prime}+(2h+1) u^{\prime}-4t u =0$. 

\par\medskip\noindent
6) P3${}^\prime$($D_7$)($\alpha, \beta, \gamma, 0$) does not have the Riccati type solution. 
P3${}^\prime$($D_7$)($0, -2, 2,0$) has an algebraic solution $y=  t^{1/3}$. 

\par\medskip\noindent
7) P3${}^\prime$($D_8$)($\alpha, \beta, 0, 0$) does not have the Riccati type solution. 
P3${}^\prime$($D_8$) ($8h,$ $-8h,0,0$) has an algebraic solution $y= -\sqrt{t} $. 

\par\medskip\noindent
8) P5($a, -a, 0, \delta$) has a rational solution $y=-1$.  
P5($ (\kappa_0+s)^2/2, -\kappa_0^2/2, -(s+1), -1/2$) has the Riccati type solutions 
$y=-tu'/ (\kappa_0+s)u$. 
Here $u$ is any solution of 
$ t^2 u^{\prime\prime}+ t( t-s-2\kappa_0+1) u^{\prime}+\kappa_0(\kappa_0+s) u 2=0$. 
If $\kappa_0=1 $, 
P5($ (s+1)^2/2, -1/2, -(s+1), -1/2$) has a rational solution $y= t/(s+1)+1$, 
which reduces to the Laguerre polynomial.

\par\medskip\noindent
9) deg-P5$( {\alpha_1^2}/2, -{\beta_1^2}/2, -2,0)$ is equivalent to 
P3($D_6$)$(4(\alpha_1-\beta_1),-4(\alpha_1+\beta_1-1), 4, \ -4)$. 
deg-P5($h^2/2,-8,-2,0$) has an algebraic solution $y=1+2\sqrt{t}/h$.
deg-P5($\alpha, 0, \gamma, 0$) has the Riccati type solutions. 

\par\medskip\noindent
10) All of classical solutions of  P1 to P5 are equivalent to the 
above solutions up to the B\"acklund transformations.
\end{theorem}

\subsection{Symmetric solutions for P1, P2, P34 and P5}
We will review symmetric solutions \cite{AVK}, \cite{KK}. 
The first, second, thirty-fourth and fourth \PA equations
\begin{eqnarray*}
&\textrm{P1}\qquad  &y^{\prime\prime}=6y^2 + t,\\
&\textrm{P2}(\alpha)\qquad  &y^{\prime\prime}= 2y^3+  t y + \alpha,\\
&\textrm{P34}(a)\qquad  & y^{\prime\prime}=\frac{(y^{\prime })^2}{2y}+2y^2 -ty -\frac{a}{2y},\\
&\textrm{P4}(\alpha,\beta)\qquad  &y^{\prime\prime}= \frac{1}{2y}{y^{\prime}}^2+\frac{3}{ 2}y^3+ 4 t y^2 
+2(t^2-\alpha)y +\frac{\beta}{y},\\
\end{eqnarray*}
have a simple symmetry:
\begin{eqnarray*}
&\textrm{P1}\qquad &y \to \zeta^3 y, \quad t \to \zeta t,  \quad (\zeta^5=1)\\
&\textrm{P2}\qquad &y \to \omega y, \quad t \to \omega^2 t, \quad  (\omega^3=1)\\
&\textrm{P34}\qquad &y \to \omega y, \quad t \to \omega t, \quad  (\omega^3=1)\\
&\textrm{P4}\qquad &y \to - y, \quad t \to - t, 
\end{eqnarray*}
There exist symmetric solutions which are invariant under the action of the cyclic group.
The symmetric solutions are studied by Kitaev \cite{AVK} for P1 and P2 and  by Kaneko 
\cite{KK} for P4.  Since these symmetric solutions exist for any parameter of the \PA equations,
they are not algebraic for generic parameters.

\begin{theorem}\label{symmetric} 1) For P1, we have two symmetric solutions
\begin{eqnarray*}
y &=& \frac16 t^3 +\frac1{336}t^8+\frac1{26208}t^{13}+\frac{95}{224550144}t^{18}+\cdots,\\
y &=& t^{-2} -\frac1{6}t^3+\frac1{264}t^{8}-\frac1{19008}t^{13}+\cdots.
\end{eqnarray*}

2) For P2($\alpha$), we have three symmetric solutions
\begin{eqnarray*}
y &=& \frac{\alpha}2 t^2 +\frac{\alpha}{40} t^5  +\frac{10\alpha^3+\alpha}{40}t^{8}+\cdots,\\
y &=&  t^{-1}-\frac{\alpha+1}{4} t^3   +\frac{(\alpha+1)(3\alpha+1)}{112}t^{5}+\cdots,\\
y &=& - t^{-1} -\frac{\alpha-1}{4} t^3  -\frac{(\alpha-1)(3\alpha-1)}{112}t^{5}+\cdots.
\end{eqnarray*}
They are equivalent to each other by the B\"acklund transformations.

2) For P34($a^2$), we have three symmetric solutions
\begin{eqnarray*}
y &=& a t  +\frac{a(2a-1)}{8} t^4  +\frac{a(2a-1)(10a-3)}{560}t^{7}+\cdots,\\
y &=& -a t  +\frac{a(2a+1)}{8} t^4 -+\frac{a(2a+1)(10a+3)}{560}t^{7}+\cdots,\\
y &=& \frac2{t^2}+\frac t2 -\frac{4a^2-9}{224} t^4 -\frac{4a^2-9}{5600}t^7+\cdots.
\end{eqnarray*}
Each solution corresponds to the symmetric solutions of P2, respectively, by 
$$y_{34}=y_2^\prime+y_2^2+\frac t2.$$

4) For P4($\alpha, -8\theta_0^2$), we have four symmetric solutions
\begin{eqnarray*}
y &=& \pm 4\theta_0\left(t -\frac{2 \alpha}{3} t^3 
     +\frac{2}{15} (\alpha^2+12\theta_0^2\pm \theta_0+1 )t^{5}+\cdots\right),\\
y &=& \pm t^{-1}+\frac{2}{3}(\pm\alpha-2) t^3 
       \mp \frac{2}{45}(-7\alpha^2\pm 16\alpha+36\theta_0^2-4 )t^{5}+\cdots.
\end{eqnarray*}
They are equivalent to each other by the B\"acklund transformations.
\end{theorem}

\par\noindent
We remark that there are no B\"acklund transformation for P1.
We may think symmetric solutions as a generalization of Umemura's classical 
solutions. We notice that symmetric solutions are not classical functions 
except for special parameters. For example, the first solution of P2($\alpha$) 
is transcendental for a generic parameter, but it is a rational solution $y=0$ 
for $\alpha=0$. It is rational if and only if $\alpha/3 $ is an integer. 
The symmetric solutions exist for 
any parameters and classical solutions exist only for special parameters.

\section{Confluent hypergeometric equation}\label{chg}
In this section we review Kummer's confluent hypergeometric equation. 
It is known that there exist two standard forms, Kummer's form and Whittaker's form  
for the confluent hypergeometric equation.  At the first we will use  Kummer's standard form.

We review   irregular singularities to fix our notations. 
For a rational function $a(x)$, we set
$${\rm ord}_{x=c}a(x)=\ ({\rm pole\ order\ of}\ a(x)\ {\rm at}\ x=c).$$
A linear differential equation 
$$\frac{d^2u}{dx^2}+p(x)\frac{d u}{dx }+q(x)=0$$
has an irregular singularity at $x=c$ if and only if
${\rm ord}_{x=c}p(x)>1$ or ${\rm ord}_{x=c}q(x)>2$. 
The Poincar\'e rank $r$ of the irregular singularity   $x=c$ is 
$$r=\max \{{\rm ord}_{x=c}p(x), 1/2 \cdot{\rm ord}_{x=c}q(x)  \}-1.$$
We think that a regular singularity has the Poincar\'e rank $r=0$.

Kummer's confluent hypergeometric equation is
\begin{equation}\label{2:chg}
x\frac{d^2u}{dx^2}+(c-x)\frac{du}{dx}-au=0, 
\end{equation}
which has a regular singularity at $x=0$ and an irregular singularity 
with the Poincar\'e rank $1$ at $x=\infty$. 
We set $x\to \varepsilon x, a \to 1/\varepsilon$ in \eqref{2:chg}  and 
take the limit $\varepsilon \to 0$. Then we get a degenerate confluent 
hypergeometric equation (the confluent hypergeometric limit equation)
\begin{equation}\label{2:dchg}
x\frac{d^2u}{dx^2}+c\frac{du}{dx}-u=0, 
\end{equation}
which has a regular singularity at $x=0$ and an irregular singularity 
with the Poincar\'e rank $1/2$ at $x=\infty$. 
The solution of \eqref{2:dchg} is
$$y=C\ {}_0F_{1}(c;x)+D\ x^{1-c}{}_0F_{1}(2-c;x).$$
\eqref{2:dchg} is reduced to \eqref{2:chg} by Kummer's second formula
$${}_0F_{1}\left(c; {x^2}/{16}\right)= e^{-x/2}\ {}_{1}F_1 (c-1/2,2c-1;x).$$
It is also related to the Bessel function 
\begin{equation}\label{bessel}
{}_0F_{1}\left(c; {x^2}/{16}\right)= \Gamma(c)(-ix/4)^{1-c}J_{c-1}(-ix/2).
\end{equation}

Later we will use  $SL$-type equations in the section\ref{pa_type}. 
$SL$-type of the confluent hypergeometric equation is called the Whittaker equation:
\begin{eqnarray}
W_{k,m}: & \dfrac{d^2u}{dx^2}= \left(\dfrac{ 1}{4}-\dfrac{k}{x}+\dfrac{m^2-\frac{1}{4}}{x^2} \right)u=0,\label{sl:chg}\\
DW_{m}:&  \dfrac{d^2u}{dx^2}= \left(\dfrac{ 1}{x}+\dfrac{m^2-\frac{1}{4}}{x^2} \right)u=0.\label{sl:dchg}
\end{eqnarray}
In $W_{k,m}$, the parameters $k, m$ correspond to 
 $k=c/2-a, \quad m= (c-1)/2$ 
in \eqref{2:chg}. In $DW_{m}$, the parameter  $m$ corresponds to 
 $ m=(c-1)/2$ in \eqref{2:dchg}.

\section{Linear equations of the \PA type}\label{iso}
We call a linear equation of the second order whose isomonodromic 
deformation gives a \PA equation 
as {\it the \PA type}. 
In this section we will list up all of linear equations of the \PA type.
Since a linear equation of the second order is equivalent to a $2\times 2$ system of 
linear equations of the first order, we use a single equation. It is easy to 
rewrite as a $2\times 2$ system. 

\subsection{Singularity type}

Linear differential equations of  the \PA type have 
the following types of singular points. 
\begin{center}
\begin{tabular}{ r|l }
P6\ &  $(0)^4$  \\
\hline
P5\ & $(0)^2(1)$ \\
\hline
P3$(D_6)$\ & $(1)^2$ or $(0)^2(1/2)$ \\
\hline
P3$(D_7)$\ & $(1)(1/2)$  \\
\hline
P3$(D_8)$\ & $(0)(2)$ \\
\hline
P4\ & $(0)(2)$ \\
\hline
P34\ &   $(0)(3/2)$ \\
\hline
P2\ & $(3)$  \\
\hline
P1\ & $(5/2)$ 
 \end{tabular} 
\end{center}
Here $(m)^k$ means $k$ singular points with the Poincar\'e rank $m$, and 
$(0)$ means a regular singularity.  The \PA third equation
has two type of singularities. $(1)^2$ is a standard one, and $(0)^2(1/2)$ is 
deg-P5. 

Two different types of linear equations are used for the isomonodromic deformation of 
the \PA second equation. One is Garnier's form \cite{Garnier}, which is the same one used by 
Okamoto \cite{Okamoto}
and Miwa-Jimbo \cite{JM2}, The second one is Flaschka-Newell's form \cite{FN}, which is 
is equivalent to the type $(0)(3/2)$   by a rational transformation \cite{Kapaev}.
It is natural that Flaschka-Newell's form as an isomonodromic deformation for P34. 
We will study Garnier's form and Flaschka-Newell's form in a succeeding paper.

\subsection{List of equations of the \PA type}\label{pa_type}

In the following, we will list up all of linear equations of the \PA type.
We use linear equations of $SL$-type:
\begin{equation}\label{SL}  
\frac{d^2 u}{d z^2}=p(z,t)u.
\end{equation}
Isomonodromic deformation for a linear equation of $SL$-type is given
by the compatibility condition of the following system \cite{Okamoto}:
\begin{equation}\label{SL:mpd}\begin{aligned} 
\frac{\partial^2 u}{\partial x^2}&=V(x,t)u,\\
\frac{\partial u}{\partial t}=& A(x,t)\frac{\partial  u}{\partial x}
                             -\frac12\frac{\partial A(x,t)}{\partial x }u.
\end{aligned}\end{equation}
We will list up $V(x,t)$ and $A(z)$ of linear equations of \PA type. We corrected 
misprints in \cite{Okamoto}.  The compatibility condition gives the \PA equation, 
which turn out a Hamiltonian system  with the Hamiltonian $K_J$ in the following list.
$(q,p)$ are canonical coordinates and $q$ satisfies the \PA equation.
\par\bigskip\noindent
{\it Type $(5/2)$}\ : the first \PA equation P1 \\
$$V(z,t)=4z^3 +2t z +2 K_{\rm I} +\frac{3}{4(z-q)^2}  -\frac{p}{z-q}$$
$$A(z)= \frac 12\cdot \frac{1}{z-q}$$
$$K_{\rm I} =\frac12 p^2 -2q^3 -t q.$$

\par\bigskip\noindent
{\it Type $(3)$}\ : the second \PA equation P2($\alpha$)\\
$$V(z,t)=z^4 +t z^2 +2\alpha z +2 K_{\rm II} +\frac{3}{4(z-q)^2}  -\frac{p}{z-q}$$
$$A(z)= \frac 12\cdot \frac{1}{z-q}$$
$$K_{\rm II}  =\frac12 p^2 -\frac12 q^4 -\frac12 t q^2 -\alpha q.$$

\par\bigskip\noindent
{\it Type $(1)(3/2)$}\ : the thirty-fourth \PA equation P34($\alpha$)\\
$$V(z,t)=\frac z2 -\frac t2 +\frac{\alpha-1}{4z^2} -\frac{ K_{\rm XXXIV}}z
 +\frac{3}{4(z-q)^2} -\frac{pq}{z(z-q)}$$
$$A(z)= - \frac{z}{z-q}$$
$$K_{\rm XXXIV}  =-qp^2 +p +\frac{q^2}2 -\frac{tq}{2} +\frac{\alpha-1}{4q}.$$

\par\bigskip\noindent
{\it Type $(1)(2)$}\ : the fourth \PA equation P4($\alpha, \beta$)\\
$$V(z,t)=\frac{{a_0}}{z^2}+\frac{K_{\rm VI}}{2z} +a_1 +\left(\frac{z+2t}4\right)^2 
+\frac{3}{4(z-q)^2}  -\frac{pq}{z(z-q)},$$
$$A(z)=   \frac{2z}{z-q},$$
$$K_{\rm VI}  =2qp^2-2p -\frac{{a_0}}{q} -2a_1 q-2q\left(\frac{q+2t}4\right)^2,$$
$$a_0= -\frac\beta 8-\frac14,\ a_1=-\frac\alpha4.$$

\par\bigskip\noindent
{\it Type $(1)^2$}\ : the third \PA equation P3${}^\prime$($\alpha, \beta, \gamma, \delta$)\\
$$V(z,t)=\frac{a_0t^2}{z^4}+\frac{a_0^\prime t }{z^3}-\frac{t K_{\rm III}^\prime}{ z^2}+
\frac{a_\infty^\prime}{ z } + a_\infty
+\frac{3}{4(z-q)^2}-\frac{pq}{z(z-q)},$$
$$A(z)=   \frac{ qz}{t(z-q)},$$ 
$$t K_{\rm III}^\prime = q^2p^2 -qp -\frac{a_0 t^2}{ q^2}
-\frac{a_0^\prime t}{ q}- {a_\infty^\prime q} -a_\infty  {q^2},$$
$$a_0= -\frac\delta{16},\ a_0^\prime=-\frac\beta8,\ 
a_\infty= \frac\gamma{16},\ a_\infty^\prime=\frac\alpha8.$$

\par\bigskip\noindent
{\it Type $(0)^2(1)$}\ : the fifth \PA equation P5($\alpha, \beta, \gamma, \delta$)\\
$$V(z,t)=\frac{a_1 t^2}{(z-1)^4}+\frac{K_{\rm V} t}{(z-1)^2 z}+\frac{a_2  t}{(z-1)^3}-\frac{p
   (q -1) q }{ z (z-1)(z-q )}+\frac{a_\infty}{(z-1)^2}+\frac{{a_0}}{z^2}+\frac{3}{4
   (z-q )^2}$$
$$A(z)=\frac{q-1}t \cdot \frac{z(z-1) }{ z-q }$$
$$ tK_{\rm V}= q(q -1)^2   \left[-\frac{a_1 t^2}{(q-1)^4}-\frac{a_2 t}{(q -1)^3}+p^2-
   \left(\frac{1}{q }+\frac{1}{q -1}\right)p -\frac{a_\infty}{(q -1)^2}-\frac{a_0}{q 
   ^2}\right] $$
$$a_0= -\frac\beta 2-\frac14,\ a_1=-\frac\delta2  ,\ a_2=-\frac\gamma2,\ a_\infty=\frac12(\alpha+\beta)-\frac34$$

\par\bigskip
If $\gamma=0$ or  $\delta=0$ for P3, the type of the linear equation is $(1)(1/2)$. 
We will take  $\delta=0$ as a standard form, which reduces to P3($D_7$)($\alpha, \beta, \gamma,0$).
If $\gamma=0$ and  $\delta=0$ for P3, the type of the linear equation is $ (1/2)^2$, which reduces to
P3($D_8$)($\alpha, \beta,0,0$). 
If $\delta=0$ for P5, the type of the linear equation is $(0)^2(1/2)$, which reduces to
deg-P5($\alpha, \beta, \gamma,0$). 
We omit P6 since we treat the cases of  irregular singularities   in this paper.

\section{Rational transform of $W_{k,m}$ and $DW_m$}\label{rational}

In this section, we will classify all of rational transformations 
of independent variables, which change $W_{k,m}$ or $DW_{m}$ into linear 
differential equations of the \PA type.

The following simple lemma is a key of classification.
\begin{lemma}
For a linear differential equation 
\begin{equation}\label{2ndnormal}
\frac{d^2u}{dx^2}=p(x)u(x),
\end{equation}
we take  a new independent variable $z$ defined by a rational transform $x=x(z)$.
If $x=c$ is a regular singularity, $z_c=x^{-1}(c)$ is also a regular singularity 
of the transformed equation. When $z_c$ is a $n$-th branched point and the distance
of the two exponent at $x=c$ is $1/n$, $z=z_c$ is an apparent singularity and 
$z=z_c$ can be reduced to a regular point by a suitable change $y=v(z)u$.

If $x=c$ is an irregular singularity of the Poincar\'e rank $m$ and 
$z_c$ is a $n$-th branched point, $z=z_c$ is an irregular singularity of the Poincar\'e 
rank $mn$.
\end{lemma}

The proof is obvious. Since we will start from the confluent hypergeometric equations, 
we take a care of $x=0, \infty$ as branched points. For $x=x(z)$, 
if  $x^{-1}(0)$ are   branch points with $\mu_1,\mu_2,\cdots$ order and 
$x^{-1}(\infty)$ are   branch points with $\nu_1,\nu_2,\cdots$ order, 
we  call  $x(z)$ as a type $(\mu_1+\mu_2+\cdots |\nu_1+\nu_2+\cdots)$.
We allow $\mu_j$ or $\nu_j$ equals to one. $\sum \mu_j= \sum\nu_j$ is 
the order of the rational function $x(z)$. 

\begin{theorem}\label{th:cover} By a rational transform $x=x(z)$, $W_{k,m}$ or $DW_m$
changes to a linear equation of the \PA type or a confluent hypergeometric equation 
if and only if one of the following cases occur. 
\par\medskip

1) Double cover 
\par\medskip
\begin{tabular}{ c|c|c|c }
 $W_{k,m}$\   & $(2|2)$     & $(0)(2)$   & P4-sym \\   	%%Kaneko
 $W_{k, 1/4}$\ & $(2|2)$     &   $(2)$    & Weber\\		%%Weber
 $W_{k, 1/4}$ \ & $(2|1+1)$   & $(1)^2 $   & D6-alg \\		%%done
 $W_{0, 1/2}$\   & $(1+1|2)$   & $(0)(2) $   & P4-Her \\		%%done
 $DW_{m}$\    & $(2|2)$ & $(0)(1) $ & Kummer \\		%%Kummer
 $DW_{m}$\    & $(1+1|2)$ & $(1)^2(2) $ & P5-rat \\		%%done
 $DW_{1/4}$\  & $(2|1+1)$ & $(1/2)^2$ & D8-alg \\		%%done
\end{tabular}

\par\medskip
2) Cubic cover
\par\medskip
\begin{tabular}{ r|c|c|c }
 $W_{k,1/3}$\   & $(3|3)$     & $(3)$   & P2-sym \\			%Kaneko
 $DW_{m}$\ & $(3|3)$      &   $(1)(3/2)$    & P34-sym \\		%Kitaev
 $DW_{1/6}$\ & $(3|3)$      &   $(3/2)$    & Airy \\		%Kitaev
 $DW_{1/4}$\ & $(2+1|3)$   &   $(0)(3/2)$    & P34-rat \\		%?
 $DW_{1/6}$\ & $(3|2+1)$   &   $(1)(1/2)$    & D7-alg  		%done
\end{tabular}
\par\medskip

3) Quartic cover
\par\medskip
\begin{tabular}{ r|c|c|c }
 $DW_{1/6}$\   & $(3+1|4)$     & \hskip 5mm  $(3) \hskip 5mm $   & P4-rat 
\end{tabular}
\par\medskip

4) Quintic cover 
\par\medskip
\begin{tabular}{ r|c|c|c }
 $DW_{1/5}$\    & \hskip 3mm $(5|5)$ \hskip 3mm {} &\hskip 3mm  $(5/2) $ \hskip 3mm &  P1-sym \\
  $DW_{1/10}$\    & \hskip 3mm $(5|5)$ \hskip 3mm {} &\hskip 3mm  $(5/2) $ \hskip 3mm &  P1-sym 				%Kitaev
\end{tabular}

5) Sextic cover 
\par\medskip
\begin{tabular}{ r|c|c|c }
 $DW_{1/6}$\    & $(3+3|6)$     &\hskip 5mm $(3)\hskip 5mm  $ & P2-rat 
\end{tabular}
\par\medskip\noindent
Here the first column is the starting linear equation. The second column is 
the type of a rational transform. The third column is the singularity type of 
of the transformed linear equation. The fourth column is the solution of 
the \PA equation.
\end{theorem}
\par\noindent
{\it Proof.}  The Poincare rank of irregular singularities of
 equations of the \PA type is at most three. Therefore we can take 
 up to a cubic cover of  $W_{k,m}$ and we can take  up to a sextic cover of  $DW_{m}$.
 In each cases, we can classify the covering map directly.
 \hfill \boxed{} 
 \par\bigskip\noindent 
The other types of confluent equations, such as the Weber (parabolic cylinder) equation
and the Airy equation, appear in the table. Therefore covering of such equations are
also included in the covering of $W_{k,m}$ or $DW_{m}$.  We consider the Bessel equation as a 
special case of $W_{0,m}$. 
We will explain each case in the next subsection.  The theorem \ref{th:cover} shows that 
it is necessary to distinguish  P2 and P34 from the isomonodromic viewpoint.

\subsection{Exact form}\label{exact}
In this section we will write down algebraic solutions obtained by rational transformations
from $W_{k,m}$ or $DW_{m}$ explicitly.
The following lemma is a key to calculate exact transformations.

\begin{lemma}\label{keylemma} For the equation
\begin{equation*} 
 \frac{d^2u}{dx^2}=Q(x)u,
\end{equation*}
we set 
$$x=x(z), \ u(x)= \sqrt{\frac{dx}{dz}} v(z).$$
Then $v$ satisfies
\begin{equation}\label{transformed}
\frac{d^2v}{dz^2} =\left( Q(x(z))(x'(z))^2- \frac12\{z,x\} \right)v.
\end{equation}
Here $\{z,x\}$ is the Schwarzian derivative
$$\{z,x\}=\frac{z^{\prime\prime\prime}}{z'}
     -\frac32 \left(\frac{z^{\prime\prime }}{z'} \right).$$
\end{lemma}
\noindent 
We will take  $Q(x)$ as \eqref{sl:chg} or \eqref{sl:dchg}. Then we calculate  $V(z,t)$ 
in the left hand side of \eqref{transformed} in each case of the theorem \ref{th:cover}. 
We should take a suitable coefficients of the rational function $x=x(z)$ 
to coincide  with the equation of the \PA type.  In the following $(q,p)$ is the 
canonical coordinates in the section \ref{pa_type}.

\subsubsection{Double cover}
1) P4-sym $(2|2)$ \\ \noindent
For $W_{k,m}$, we set $x=z^2/4$.  Then  
$$V(z,t)= -k+ \frac{16m^2-1}{4 z^2}+\frac {z^2}{16},$$
which is the case $t=0$ for the symmetric solution $ q(t)=2(4m+1)t+O(t^3)$ 
of P4$(4k,-2(4m+1)^2)$.

\par\bigskip\noindent
2) Weber $(2|2)$ \\ \noindent
This is a special case of P4-sym.
For $W_{k,1/4}$, we set $x=z^2/2$.  Then  
$$V(z,t)= \frac{z^2}{4}  -2k,$$
which is the parabolic cylinder equation for $D_{2k-1/2}(z)$. 
This shows the  well-known formula
\begin{equation}\label{weber}
D_{2k-1/2}(z)= 2^k z^{-1/2}W_{k,-1/4}\left(\frac{z^2}2 \right).
\end{equation}

\par\bigskip\noindent
3) $D_6$-alg $(2|1+1)$ \\ \noindent
For $W_{k,1/4}$, we set $x=(z-\sqrt{t})^2/z$.  Then  
\begin{eqnarray*}
V(z,t)= \frac{1}{4} + \frac{t^2}{4z^2}-\frac{kt}{z^3}
-\frac{8t+32k\sqrt{t}+3}{16 z^2}
-\frac{k}{z }+ \frac{3}{ 4(z+\sqrt{t})^2}-\frac{3}{ 4z(z+\sqrt{t})},
\end{eqnarray*}
which gives the algebraic solution $ q(t)=-\sqrt{t}$ of P3$(-8k,8k,4,-4)$.

\par\bigskip\noindent
4) P4-Her $(1+1|2)$ \\ \noindent
For  $W_{0, 1/2}$, we set $x =z (z+4t)/4$. Then  
$$V(z,t)= \left(\frac{z+2t}{4} \right)^2+\frac{4}{ 3(z+2t)},$$
which gives the rational solution $ q(t)=-2t$ of P4($0,-2$). 

\par\bigskip\noindent
5) Kummer's second formula $(2|2)$ \\ \noindent
For $DW_{m}$, we set $x=z^2/16$.  Then  
$$V(z,t)= \frac{1}{4} + \frac{16m^2-1}{4 z^2},$$
which is $W_{0, 2m}$. This is Kummer's second formula.

\par\bigskip\noindent
6) P5-rat $(1+1|2)$  \\ \noindent
For $DW_{m}$, we set $x= h t^2 z  /4(z-1)^2$.  Then  
$$V(z,t)= \frac{h t^2}{(z-1)^4} + \frac{16m^2-1+h t^2}{4 z(z-1)^2}
-\frac3{4z} + \frac{4m^2-1}{4z^2} + \frac{3(z+2)}{4(z+1)^2},$$
which gives the rational solution $ q(t)=-1$ of  P5($2m^2, -2m^2, 0, -2h$).

\par\bigskip\noindent
7) $D_8$-alg $(2|1+1)$ \\ \noindent
For $DW_{1/4}$, we set $x=h (z-\sqrt{t})^2/z$.  Then  
$$V(z,t)= \frac {h t}{z^3} +\frac {32 h \sqrt{t}-3}{16z^2}
+\frac{h}z+ \frac{3}{4(z+\sqrt{t})^2}-\frac{3}{4z(z+\sqrt{t})},$$
which gives the algebraic solution $ q(t)=-\sqrt{t}$ of P3$(8h,-8h,0,0)$.

\subsubsection{Cubic cover}
1) P2-sym $(3|3)$\\ \noindent
For $W_{k,1/3}$, we set $x=2z^3/3$.  Then  
$$V(z,t)= z^4 -6kz +\frac{3}{4z^2},$$
which is the symmetric solution $ q(0)=0, p(0)=0$ of P2$(-3k)$.

\par\bigskip\noindent
2) P34-sym $(3|3)$\\ \noindent
For $DW_{m}$, we set $x= z^3/18 $.  Then  
$$V(z,t)=\frac z 2 +\frac{36m^2-1}{4z^2},$$
which is the symmetric solution $ q(0)=0$ of P34$(3 (12m^2 -1))$.  
We have two symmetric solutions with $ q(0)=0$, and both give the same 
equation when $t=0$.

\par\bigskip\noindent
3) Airy $(3|3)$\\ \noindent 
This is a special case of P34-sym.
For $DW_{1/6}$, we set $x=  z^3/9 $.  Then  
$$V(z,t)= z,$$
which is  the Airy equation. It is know that 
\begin{equation}\label{airy}
{\rm Ai}(x)=\frac1{ 3^{2/3} \Gamma(\frac23) } {}_0F_1\left(\frac23; \frac{z^3}9  \right)
-\frac{x}{3^{1/3}\Gamma(\frac13)}{}_0F_1 \left(\frac43; \frac{z^3}9  \right).
\end{equation}

\par\bigskip\noindent
4) P34-rat $(2+1|3)$\\ \noindent 
For $DW_{1/4}$, we set $x=  z(z-3t/2)^3/18 $.  Then  %% h=1/18 k=-1/2
$$V(z,t)=  \frac z 2-\frac  t2 +\frac{t^3/4+1}{2tz}-\frac{3}{16z^2} 
          +\frac{3}{4(z-t/2)^2}-\frac{1}{2t(z-t/2)},$$
which gives the rational solution $ q(t)=t/2$ of P34$( 1/4)$.

\par\bigskip\noindent
5) $D_7$-alg $(3|2+1)$\\ \noindent
For $DW_{1/6}$, we set $x =\left(z+2 t^{1/3} \right)^3/32z$. Then  
$$V(z,t)= \frac{ t}{4z^3}-\frac{16+27 t^{2/3}}{72z^2}
+\frac{2}{ 3 t^{1/3} z }+ \frac18   +\frac{3}{4 (  z-  t^{1/3})^2}
- \frac{2}{3t^{1/3}(  z-  t^{1/3}) },$$
which is the algebraic solution $ q(t)=  t^{1/3}  $ of P3($0, -2, 2,0$).

\subsubsection{Quartic cover}
If the covering degree is more than three, we start only from $DW_{m}$. 
If  a quartic covering map $x=x(z)$ splits to $x= (\bar{x}(z))^2$, 
the covering through  $W_{0,2m}$, which is Kummer's case 5.1.4. 
Such case occurs in case that all of $\mu_j, \nu_j$ are even.
$(4|4), (4|2+2), (2+2|4)$ corresponds to $(2|2), (2|1+1), (1+1|2)$, respectively.

\par\bigskip\noindent
1) P4-rat $(3+1|4)$ \\ \noindent
For $DW_{1/6}$, we set  $   x=\frac{1}{256} z \left(z+  8 t/3 \right)^3$.
Then 
$$V(z,t)= -\frac{2}{9 z^2}+\frac{2t^3}{27 z} +\left(\frac{z+2t}4\right)^2
+\frac{27}{4(3z+2t)^2}-\frac1{z(3z+2t)}$$ 
which gives the rational solution $q(t)=-2t/3$ of P4$(0, -2/9)$.

\subsubsection{Quintic cover}
1) P1-sym $(3+1|4)$ \\ \noindent
For $DW_{1/5}$, we set $x =4z^5/25$. Then  
$$V(z,t)= \frac{3}{4 z^2}+4 z^3,$$
which is the symmetric solution $ y(0)=0, y'(0)=0, t=0$ of P1.
\par\bigskip \noindent
For $DW_{1/10}$, we set $x =4z^5/25$. Then  
$$V(z,t)=  4 z^3,$$
which is the symmetric solution $ y(t)=1/t^2+\cdots$ and $t=0$ of P1. 
If we substitute $ y(t)=1/t^2+\cdots$ in $Q(z,t)$, $Q(z,t)$ may have a pole 
$t=0$ but this pole is apparent. 

\subsubsection{Sextic cover}
As the same as the quartic cover, we omit the case $(6|6), (4+2|6)$. 
\par\bigskip\noindent
1) P2-rat $(3+3|6)$ \\ \noindent
For $DW_{1/6}$, we set $x =\frac1{36} (z^2+t)^3$. Then  
$$V(z,t)= \frac{3}{4 z^2}+tz^2+ z^4,$$
which is the rational solution $q=0$ of P2$(0)$.

\subsection{R.~Fuchs' Problem}

Compared with the theorem \ref{classical} and the theorem \ref{th:cover}, 
we obtain all of algebraic solutions and symmetric solutions except 
algebraic solutions of deg-P5 and the Laguerre solutions of  P5. 
These two solutions are not obtained from $W_{k,m}$ or $DW_{m}$ by rational transformations, 
but it is  obtained by exponential type transformations. 

\subsubsection{Split case}
1) For deg-P5-alg, we start from
$$\frac{d^2 u}{dx^2} =\frac {h^2-1}{4 x^2}u.$$
We set $x=e^{4\sqrt{tz/(z-1)}/h }(\sqrt{z}+\sqrt{z-1} )/ (\sqrt{z}-\sqrt{z-1} )$. Then 
\begin{eqnarray*}
V(z,t)=\frac{t}{ (z-1)^3}-\frac{4h^2-13 }{16(z-1)^2}-\frac3{16z^2}
        -\frac{2(h+2\sqrt{t})^2-5}{8z(z-1)^2}\\
   +\frac{3}{4(z-2\sqrt{t}/h-1)^2}-\frac{3+8\sqrt{t}/h}{4z(z-1)(z-1-2\sqrt{t}/h)},
\end{eqnarray*} 
which gives an algebraic solution $y=1+2\sqrt{t}/h$ for P5($h^2/2,-8,-2,0$).

\par\bigskip\noindent
2)  For P5-Lag, we start from
$$\frac{d^2 u}{dx^2} =\frac {h^2-1}{4 x^2}u.$$
We set $x=e^{t/(h(z-1))}(z-1) $. Then 
$$V(z,t)=\frac{t^2}{4(z-1)^4}-\frac{ht}{2(z-1)^3}+\frac{h^2/4-1}{(z-1)^2}
         -\frac{3}{4(z-t/h-1)^2}-\frac{ht}{(z-1)^2(z-t/h-1)},$$
which gives a rational solution $y= t/h+1$ for P5($h^2/2,-1/2,-h,-1/2$).

In these two cases, the monodromy group of  
\begin{equation}\label{split}
\frac{d^2u}{dz^2}= V(z,t)u(z)
\end{equation} 
is diagonal.   Therefore we cannot reduce them to $W_{k,m}$ or $DW_{m}$.

\subsubsection{Summary}
We summarize our result:
\begin{theorem}\label{final} By  any rational transformation from confluent hypergeometric 
equations $W_{k,m}$ or $DW_{m}$ to equations of the \PA type,  we get algebraic solutions or symmetric 
solutions of P1, P2, P34, P3, P4, P5 except deg-P5-alg and P5-Lag. 
Conversely, any algebraic solution except deg-P5-alg and P5-Lag or symmetric solution  
from P1 to P5 is obtained by  a rational transformation of $W_{k,m}$ or $DW_{m}$. 
For deg-P5-alg and P5-Lag, we can reduce them to  differential equations with 
 constant coefficients by exponential type transformations.
\end{theorem}

\par\noindent
In this sense, R.~Fuchs' problem is true for P1 to P5.  We obtain rational 
solutions and symmetric solutions of P2 and  P34 in different way. 
In \cite{AVK}, Kitaev studied a symmetric solution of P2, but his result is 
on   a symmetric solution of P34 from our viewpoint since he used  Flaschka-Newell's form. 
Symmetric solutions of P2  is studied in \cite{KK2} by using   Miwa-Jimbo's form.

The authors do not know R.~Fuchs' problem is true or not in the case the sixth \PA equations. 
Recent work by A.~V.~Kitaev \cite{Kitaev:05}  give partial, but affirmative 
answers to R.~Fuchs' problem.
We do not know any negative examples of R.~Fuchs' problem.

\end{document}